\documentclass{article}
\usepackage{spconf}
\usepackage{verbatim}
\usepackage{color}
\usepackage{amsthm}
\usepackage{amsmath}
\usepackage{amssymb}
\usepackage{graphicx}
\usepackage{mathrsfs}
\usepackage{mdwlist}
\usepackage{amsthm}
\usepackage[ruled]{algorithm2e}
\usepackage{url}
\usepackage{cite}
\usepackage{epstopdf}
\usepackage[latin9]{inputenc}
\usepackage[active]{srcltx}
\usepackage{graphicx}
\usepackage{caption}
\usepackage{subcaption}
\newtheorem{theorem}{Theorem}

\newtheorem{definition}{Definition}

 \newcommand{\norm}[1]{\left\lVert#1\right\rVert}
\title{\vspace{-1.4cm}Distributed Nonconvex Optimization for Sparse Representation\vspace{-.2cm}}
\name{ Ying Sun and  Gesualdo Scutari\thanks{Sun, and Scutari are  with the School of Industrial Engineering, Purdue University, West-Lafayette, IN, USA; emails: \texttt{<sun578,gscutari>@purdue.edu}. This work  was supported by  the USA NSF  Grants CIF
1564044, CCF 1632599, and CAREER Award 1555850, and the ONR N00014-16-1-2244. } \vspace{-.6cm}} 
\address{}
\begin{document}
\ninept
\maketitle
\begin{abstract}\vspace{-0.1cm}
We consider a   non-convex constrained Lagrangian formulation of a fundamental  bi-criteria optimization
problem for variable selection in statistical learning; the two criteria are a smooth (possibly) \emph{nonconvex} loss function, measuring the
fitness of the  model to data, and the latter  function is a \emph{difference-of-convex} (DC)  regularization,  employed to promote some extra structure on the solution,  like sparsity.  This general class of \emph{nonconvex}  problems
arises in many big-data applications, from statistical  machine learning   to physical sciences and engineering. We
develop the first  unified \emph{distributed} algorithmic framework for these problems and establish its asymptotic convergence to d-stationary solutions. Two key features of the method are: i) it can be implemented on \emph{arbitrary}  networks (digraphs) with (possibly)  time-varying connectivity; and ii) it  does not require the restrictive assumption that the (sub)gradient of the objective function is  bounded, which enlarges significantly the class of statistical learning problems that can be solved with convergence guarantees. \vspace{-0.1cm} \end{abstract}
\begin{keywords}
\! Distributed statistical learning, nonconvex optimization, sparse representation, time-varying network.
\vspace{-0.2cm}
\end{keywords}
\section{Introduction}\vspace{-0.2cm}
\label{sec:intro}
Sparse representation  \cite{Hastie_book15} is a fundamental methodology of data science in solving a broad range of
problems from statistical  machine learning  to physical sciences and engineering. 
 Significant advances have been
made in the last decade on constructing intrinsically low-dimensional solutions in high-dimensional
problems via convex programming \cite{bach2012optimization,zhang2015survey,Hastie_book15}, due to  its favorable theoretical guarantees and many efficient solution methods.
 Yet there
is increasing evidence supporting the use of non-convex formulations to enhance the realism of the
models and improve their generalizations \cite{YinLouHeXin15,LouYinXin16,ZhangXin16}. For instance, in compressed sensing, it is well documented that \emph{nonconvex} surrogates of the $\ell_0$ norm (e.g., the difference of $\ell_1$ and $\ell_2$ \cite{YinLouHeXin15},  the SCAD \cite{fan2001variable}, the ``transformed" $\ell_1$  penalty \cite{ZhangXin16}) outperform the renowned $\ell_1$ norm.  Motivated by this new line of works, in this paper, we  formulate the problem of learning a sparse parameter $\mathbf{x}$ of a statistical model from a training data set $\mathcal{D}$ as the following  general \emph{nonconvex}   constrained Lagrangian-based bi-criteria   minimization\vspace{-0.2cm}
	\begin{equation}
	\begin{aligned}
	\!\!\!\min_{\mathbf{x}\in \mathcal{K}} U\left(\bf x;\mathcal{D}\right)\triangleq \underset{F\left(\mathbf{x};\mathcal{D}\right)}{\underbrace{\sum_{i=1}^I f_i\left({\bf x};\mathcal{D}_i\right)}} + \lambda\cdot  \underset{G\left(\mathbf{x}\right)}{\underbrace{\left(G^+\left(\bf x\right)-G^-\left(\bf x\right)\right)}},
	\end{aligned}\tag{P}\label{eq: P} \vspace{-0.2cm}
	\end{equation}
where  $f_i$  is a smooth (possibly) nonconvex function measuring the fitness of the learning model to (a portion of ) the data set $\mathcal{D}_i\subseteq\mathcal{D}$;   $G$ is a penalty function, having a DC structure with $G^+$  and $G^-$ being   (possibly) nonsmooth   and  smooth, respectively;   $\lambda > 0$ is  a  parameter balancing the model fitness and  sparsity of the solution; and $\mathcal{K}\subseteq \mathbb{R}^m$ is a closed, convex set (not necessarily bounded).
 
 Problem \eqref{eq: P} is very general and encompasses a variety of \emph{convex} and \emph{nonconvex}    statistical learning formulations, including    least squares, logistic regression,    maximum likelihood  estimation,      principal component analysis,  canonical component analysis, and  low-rank approximation \cite{friedman2001elements}, just to name a few.  Furthermore, the DC structure of the penalty function $G$ allows  to accomodate in an unified fashon   either convex or nonconvex sparsity-inducing surrogates  of the $\ell_0$ cardinality function; examples are  the $\ell_p$  ($p\geq 1$), $\ell_{1,2}$ norm, the total variation penalty \cite{tibshirani1996regression,yuan2006model,rudin1992nonlinear},   the SCAD \cite{fan2001variable} function, the logarithmic \cite{weston2003use}/exponential \cite{bradley1998feature} functions, and the $\ell_p$ norm with $0<p<1$ \cite{fu1998penalized} (cf.\,Sec.\,\ref{problem_statement} for details). \\\indent
Common to the majority of the aforementioned learning tasks is the prohibitively large size of the data set $\mathcal{D}$.\,Furthermore, in several
scenarios, e.g.,  cloud, sensor, or cluster-computer networks, data $\mathcal{D}_i$'s are not centrally available but  spread  (stored) over a network, and collecting them  can be challenging or even impossible, owing
to the size of the network and volume of data, time-varying connectivity,
energy constraints, and/or privacy issues. All in all, the aforementioned reasons motivate the design of  reduced-complexity \emph{decentralized} algorithms.  This paper addresses this task. Specifically, we consider a network of $I$ agents (nodes), each of them owing a portion $\mathcal{D}_i$ of the data set $\mathcal{D}$. The network is modeled as arbitrary (possibly) time-varying digraph.  Designing   distributed solution methods for the class of problems \eqref{eq: P} in the aforementioned setting  poses several challenges, namely:    
 i) $U$ is \emph{nonconvex}, \emph{nonsmooth}, and \emph{nonseparable}; moreover, each agent $i$ knows only its own function $f_i$ [data $\mathcal{D}_j$, $j\neq i$, are not available to agent $i$]; ii) the network digraph is \emph{time-varying}, with \emph{no specific} structure; and  iii) the (sub)gradient of $U$  may not be bounded. Current works cannot address all the above challenges, as briefly
documented next.
\\\indent Most of the literature  on distributed   optimization deals with
 \emph{convex, unconstrained} optimization problems  \cite{nedic2009distributed,gharesifard2014distributed,nedich2016achieving,xu2015augmented} over \emph{undirected, static} graphs  {\cite{shi2015extra,mota2013d,duchi2012dual}}. 
 The nonconvex case has been recently studied in    \cite{bianchi2011convergence,Lorenzo2016NEXT, tatarenko2015non,sun2016distributed}. All these works however  require that the (sub)gradient of the objective function is bounded, an assumption that is not satisfied  by many formulations (e.g., least squares).  Furthermore,  \cite{tatarenko2015non} consider only  unconstrained problems, and \cite{Lorenzo2016NEXT} is applicable only to  specific network topologies (e.g., digraphs 
 with  a doubly stochastic adjacency matrix), which  limits its practical   applicability \cite{gharesifard2010does}. 
   \\\indent In this paper we address   all challenges  i)-iii) and propose the first distributed algorithmic framework for the general class of problems   (\ref{eq: P}). To cope with  i) and ii) we introduce  a  general
convexification-decomposition technique that hinges on our recent
 SCA methods \cite{facchinei2015parallel,ScuFacSonPalPan2014,daneshmand2015hybrid},  coupled with a gradient tracking mechanism, instrumental to locally estimate the missing global information.  After updating their local copy of the common variables $\mathbf{x}$,  all  agents communicate some    information to their neighbors.  
This is done using a broadcast protocol  that   requires neither a specific network topology nor  the use  of   double-stochastic consensus matrices to work [addressing thus challenge ii)]; only column stochasticity is needed.     Asymptotic convergence to d-stationary solutions
of  (\ref{eq: P}) is established, without requiring any boundedness  of   the (sub)gradient of $U$ [challenge iii)].
 Preliminary numerical
results, show that the proposed scheme compare favorably with 
state-of-the-art algorithms.\vspace{-0.2cm}

\section{Distributed Learning Model}\label{problem_statement}\vspace{-0.2cm}
We study Problem \eqref{eq: P} under the following blanket assumptions. \\
\noindent\textbf{Assumption A (Problem Setup)}\vspace{-0.1cm} 
\begin{enumerate*}
	\item[\textbf{(A1)}] The set $\mathcal{K}\neq\emptyset$ is closed and convex;
	\item[\textbf{(A2)}] Each $f_i:\mathcal{O}\rightarrow \mathbb{R}$ is    $C^1$, with Lipschitz continuous gradient $\nabla f_i$ on $\mathcal{K}$, where $\mathcal{O}\supseteq \mathcal{K}$ is an open set;
	\item[\textbf{(A3)}]   $G^+:\mathcal{K}\rightarrow \mathbb{R}$ is convex (possibly) nonsmooth and $G^-:\mathcal{O}\rightarrow \mathbb{R}$ is convex  $C^1$, with   gradient  $\nabla G_i^-$ Lipschitz  on $\mathcal{K}$;
	\item[\textbf{(A4)}] Problem \eqref{eq: P} has a solution. 
\end{enumerate*}
 
 Assumptions above are quite general and   satisfied by several loss and penalty functions, proposed in the literature. For instance,  (nonconvex) quadratic, Huber, and logistic loss functions fall under  A2. A3 is satisfied by  (nonsmooth)  convex  functions  and the majority of sparsity-inducing nonconvex surrogates of the $\ell_0$ norm proposed  up to date;  Table\,1 summarizes the  majority of the latter functions. One can see that  all   functions $G$ therein are separable, $G\left(\mathbf{x}\right)\triangleq \sum_{j=1}^{m}g\left(x_j\right)$, with  
  $g:\mathbb{R}\rightarrow \mathbb{R}$ having  the following DC form \cite{le2015dc}\vspace{-0.1cm}
 \begin{equation}
 g\left(x\right)=\underset{\triangleq g^{+}\left(x\right)}{\underbrace{\eta\left(\theta\right)\left\vert x\right\vert }}-\underset{\triangleq g^{-}\left(x\right)}{\underbrace{\left(\eta\left(\theta\right)\left\vert x\right\vert -g\left(x\right)\right)}},\label{eq: dc decompose}\vspace{-0.2cm}
 \end{equation}
 where $\eta\left(\theta\right)$ is a given function,  whose  expression depends on the surrogate $g$ under consideration, see Table\,2.   It can be shown that for all the functions in Table\,1, $g^-\left(x\right)$  has Lipschitz continuous gradient \cite{le2015dc} (the closed form is given in  Table 2), implying that A3 is satisfied.  We conclude this list of examples satisfying A1-A4, with two  concrete sparse representation problems.  
 
 \noindent\textbf{Example \#1 (Sparse Linear Regression):} Consider the problem of retrieving a sparse signal $\mathbf{x}$ from observations $\{\mathbf{b}_i\}_{i=1}^{I}$, where each $\mathbf{b}_i = \mathbf{A}_i\mathbf{x}$ is a linear measurement of the signal. The problem reads\vspace{-0.2cm}
 \begin{equation}
 \begin{aligned}
 &\underset{\mathbf{x}}{\textrm{min}}& & \sum_{i=1}^{I}\norm{\mathbf{b}_i - \mathbf{A}_i\mathbf{x}}^2 + \lambda G\left(\mathbf{x}\right),
 \end{aligned}\label{p: sparse regression}\vspace{-0.2cm}
 \end{equation}
where $G$ can be  any of the penalty functions discussed above. For instance,   if  $G$ is the $\ell_2$ and $\ell_1$ norm, \eqref{p: sparse regression} reduces to the  ridge and LASSO regression, respectively. Problem \eqref{p: sparse regression} is clearly an instance of \eqref{eq: P} with    $\mathcal{D}_i \triangleq  \{\left(\mathbf{A}_i,\mathbf{b}_i\right)\}$ and $f_i\left(\mathbf{x},\mathcal{D}_i\right) \triangleq \norm{\mathbf{b}_i - \mathbf{A}_i\mathbf{x}}^2$.  
 
 \noindent\textbf{Example \#2 (Sparse PCA):} Consider  finding the sparse principal component of a data set given by the rows of matrices $\mathbf{D}_i$'s, which leads to\vspace{-0.2cm}
 \begin{equation}
 \begin{aligned}
 &\underset{\norm{\mathbf{x}}_2 \leq 1}{\textrm{max}}& & \sum_{i=1}^{I}\norm{\mathbf{D}_i\mathbf{x}}^2 - \lambda G\left(\mathbf{x}\right),\\
 \end{aligned} \label{p: sPCA}\vspace{-0.2cm}
 \end{equation}
 where $G$ is some regularizer satisfying A3. Clearly, \eqref{p: sPCA} is a (nonconvex) instance of  Problem \eqref{eq: P}, with  $\mathcal{D}_i \triangleq \{\mathbf{D}_i\}$ and $f_i \triangleq -\norm{\mathbf{D}_i\mathbf{x}}^2$.\smallskip
 \begin{table}[t] 
 		\begin{tabular}{ c | l   }
 			\hline			
 			Penalty function & \multicolumn{1}{|c}{Expression}     \\ \hline \vspace{-0.2cm}\\ 			
 			Exp \cite{bradley1998feature} & $g_{\text{exp}}(x)= 1-e^{-\theta |x|}$   \\
 			
 			$\ell_p(0<p<1)$ \cite{fu1998penalized} & $g_{\ell_p^{+}}(x)= (|x|+\epsilon)^{1/\theta}$,  \\
 			
 			$\ell_p(p<0)$ \cite{rao1999affine} & $g_{\ell_p^{-}}(x)= 1-(\theta |x|+1)^{p}$ \\

 			SCAD \cite{fan2001variable}& $g_{\text{scad}}(x)\!= \!\begin{cases}
 			
 			\frac{2 \theta}{a+1} |x|,  & 0 \leq |x|\leq \frac{1}{\theta}  \\
 			
 			\frac{-\theta^2 |x|^2+2 a \theta |x|-1}{a^2-1},& \frac{1}{\theta} < |x|\leq \frac{a}{\theta}  \\
 			
 			1, & |x|> \frac{a}{\theta}  \end{cases}$  \\
 			
 			Log \cite{weston2003use} & $g_{\log}(x)=\frac{\log(1 +\theta |x|)}{\log(1 +\theta)}$ \\
 			
 			
 			
 			\hline
 		\end{tabular}
 		\caption{Examples of DC penalty functions satisfying A3 [cf.\,(\ref{eq: dc decompose})]
 		}\label{tab: ncvx regularizer}  \vspace{-0.4cm}
 	\end{table}
 \noindent\textbf{Network Topology.} Time is slotted and, at each time-slot $n$, the communication network of agents is  modeled as a  (possibly) time-varying digraph  $\mathcal{G}\left[n\right]=\left(\mathcal{V},\mathcal{E}\left[n\right]\right))$, where the set of vertices $\mathcal{V}=\{1,\ldots,I\}$ represents the set of  $I$ agents, and the set of edges $\mathcal{E}\left[n\right]$ represents the agents' communication links. The  in-neighborhood of agent $i$ at time $n$ (including node $i$) is defined as  $\mathcal{N}_i^{\rm in}[n]=\{j|(j,i)\in\mathcal{E}[n]\}\cup\{i\}$ whereas its out-neighbor is defined as $\mathcal{N}_i^\textrm{out}\left[n\right]=\{j|\left(i,j\right)\in\mathcal{E}\left[n\right]\}\cup\{i\}$. The out-degree of agent $i$ is defined as $d_i\left[n\right] \triangleq  \left|\mathcal{N}_i^\textrm{out}\left[n\right]\right|$.  To let information propagate over the network,  we assume that  the graph sequence $\left(\mathcal{G}\left[n\right]\right)_{n\in\mathbb{N}}$ possesses some ``long-term'' connectivity property, as formalized next.\smallskip\\
\noindent \textbf{Assumption B (On the graph connectivity).}  
The graph sequence $\{\mathcal{G}[n]\}_{n\in\mathbb{N}}$ is $B$-strongly connected, i.e., there exists an integer $B > 0$  (possibly unknown to the agents) such that the graph with edge set $\cup_{t=kB}^{(k+1)B-1} \mathcal{E}[t]$ is strongly connected,  for all $k\geq0$.
\smallskip

As a non-convex optimization problem, globally optimal solutions of \eqref{eq: P}  are in general not possible to be computed. Thus, one has to settle for computing a ``stationary'' solution in practice.  Among all the   stationarity concepts,   arguably,  a d(irectional)-stationary solution 
is  the sharpest kind of stationarity for the class of convex constrained nonconvex nonsmooth problem \eqref{eq: P}; see, e.g., \cite{pang2015computing}.  \vspace{-.1cm}
\begin{definition}[d-stationarity]
	A point $\mathbf{x}^\ast\in \mathcal{K}$ is a d-stationary solution  of \eqref{eq: P}  if $U'\left(\mathbf{x}^\ast;\mathbf{x}-\mathbf{x}^\ast\right) \geq 0$, $\forall \mathbf{x}\in \mathcal{K}$, where $U'\left(\mathbf{x}^\ast;\mathbf{x}-\mathbf{x}^\ast\right) $ is the directional derivative of $U$ at $\mathbf{x}^\ast$ along the direction $\mathbf{x}-\mathbf{x}^\ast$.\vspace{-.2cm}
\end{definition}
Quite interestingly, such nonzero d-stationary solutions have been proved to   possess some sparsity (and even minimizing) property, under a set of assumptions, including a specific choices of  $F$, $G$, and $\lambda$ in \eqref{eq: P}; we refer to \cite{pang2016PartI} for details. Motivated by these results,   our goal is then to  devise   algorithms   converging to   d-stationary solutions of Problem \eqref{eq: P}, in the above distributed setting.\vspace{-.2cm}
 \begin{table}[t]
	\begin{tabular}{ c | c|  l}
		\hline	 
		$g$ & $\eta(\theta)$ & \multicolumn{1}{c}{$\nabla g^{-}_{\theta}(x)$ }   \\\hline   \vspace{-0.3cm}	\\
		$g_{\textrm{exp}}$ & $\theta$ & $\text{sign}(x)\cdot\theta\cdot(1-e^{-\theta |x|})$  \\		
		$g_{\ell_p^{+}}$ &  $\frac{1}{\theta}\epsilon^{1/\theta-1}$ & $\frac{1}{\theta}\,\text{sign}(x)\cdot[\epsilon^{\frac{1}{\theta}-1}-(|x|+\epsilon)^{\frac{1}{\theta}-1}]$ \\
		
		$g_{\ell_p^{-}}$ & $-p \cdot\theta$  & $-\text{sign}(x)\cdot p\cdot \theta\cdot [1-(1+\theta |x|)^{p-1}]$ \\

		$g_{\text{scad}}$ & $\frac{2 \theta}{a+1}$  & $\begin{cases}
				0,  & |x|\leq \frac{1}{\theta}  \\
				
				\text{sign}(x)\cdot \frac{2 \theta(\theta |x|-1)}{a^2-1}, & \frac{1}{\theta} < |x|\leq \frac{a}{\theta}  \\
				
				\text{sign}(x)\cdot\frac{2 \theta}{a+1}, & \text{otherwise}
		\end{cases} $\\
		
		$g_{\log}$ &    $\frac{\theta}{\log(1+\theta)}$ & $\text{sign}(x)\cdot \frac{\theta^2 |x|}{\log(1+\theta)(1+\theta |x|)}$ \\
		
		
		
		
		\hline
	\end{tabular}
	 \caption{ Explicit expression of $\eta(\theta)$  and   $\nabla g^{-}(x)$ [cf.\,(\ref{eq: dc decompose})]   }\label{tab: gradient}
\vspace{-0.4cm}\end{table}
\section{Algorithm Design}\vspace{-.2cm}

We start introducing an informal description of the algorithm that sheds light on the main ideas behind the proposed framework.

Each agent $i$ maintains  a  copy of the common optimization variable $\mathbf{x}$, denoted by $\mathbf{x}_{(i)}$, which needs to be updated locally at each iteration so that asymptotically 1) $\mathbf{x}_{(i)}$ reaches a d-stationary point of Problem \eqref{eq: P}; and 2) all $\mathbf{x}_{(i)}$'s are consensual, i.e., $\mathbf{x}_{(i)} = \mathbf{x}_{(j)}, \forall i\neq j$. To do so, the proposed algorithm framework combines SCA
techniques (Step 1 below) with a consensus-like step implementing
a novel broadcast protocol (Step 2), as described next.

\noindent\textbf{Step 1: Local SCA.} 
At iteration $n$,   agent $i$ should solve  \eqref{eq: P}. However, $F - G^-$  is nonconvex and  $\sum_{j\neq i}f_j$ in $F$ is unknown. To cope with these issues,   agent $i$ solves instead    an approximation of  \eqref{eq: P}  wherein $F - G^-$   is replaced by the   strongly convex function $\widetilde{F}_i$:\vspace{-0.2cm}
\begin{equation}
\begin{aligned}
	\!\!\! \widetilde{F}_i\left(\mathbf{x}_{(i)};\mathbf{x}_{(i)}^n\right) \triangleq & \,\widetilde{f_i}\left(\mathbf{x}_i;\mathbf{x}_{(i)}^n\right)  - \nabla G^-\left(\mathbf{x}_{(i)}^n\right) ^\top \left(\mathbf{x}_i - \mathbf{x}_{(i)}^n\right)\\
	  & 
	 + \widetilde{\boldsymbol{\pi}}_i  ^{n\top} \left(\mathbf{x}_i - \mathbf{x}_{(i)}^n\right),
	\end{aligned}\vspace{-0.2cm}
\end{equation}
where   $\widetilde{f}_i: \mathcal{K}\to \mathbb{R}$ should be regarded as a (simple) strongly
convex approximation  of $f_i$ at the current iterate $\mathbf{x}_{(i)}^n$  that preserves the first order properties of $f_i$  (see Assumption C below); the second term is the linearization of the concave smooth function $-G^-$; and the last term accounts for the lack of knowledge of   $\sum_{j\neq i}f_j$:   $\widetilde{\boldsymbol{\pi}}_i^n$ aims to  track  the gradient  of $\sum_{j\neq i}f_j$. In Step 2 we will show how to update  $\widetilde{\boldsymbol{\pi}}_i^n$ so that  $\|\widetilde{\boldsymbol{\pi}}_i^n - \sum_{j\neq i}\nabla f_j\left(\mathbf{x}_{(i)}^n\right)\|\underset{n\rightarrow\infty}{\longrightarrow}0$ while using only \emph{local} information. We require $\widetilde{f}_i$
to satisfy the following mild assumptions ( $\nabla \widetilde{f}_i$ is the partial gradient of $\widetilde{f}_i$ w.r.t. the first argument). 
  
\noindent \textbf{Assumption C (On the surrogate function $\widetilde{f}_i$).}\vspace{-0.1cm}
\begin{enumerate*}
	\item[\textbf{(C1)}] $\nabla \widetilde{f}_i\left(\mathbf{x};\mathbf{x}\right)=\nabla f_i\left(\mathbf{x}\right)$, for all $\mathbf{x}\in \mathcal{K}$;
	\item[\textbf{(C2)}] $\widetilde{f}_i\left(\bullet;\mathbf{y}\right)$ is uniformly strongly convex  on $\mathcal{K}$;
	\item[\textbf{(C3)}] $\nabla \widetilde{f}_{i}\left(\mathbf{x};\bullet\right)$ is uniformly Lipschitz continuous on $\mathcal{K}$.\vspace{-0.1cm}
\end{enumerate*}

A wide array of choices for   $ \widetilde{f}_i$ satisfying Assumption\,C can be found in \cite{facchinei2015parallel}, see Sec.\,\ref{subsec:discussion} for some significant examples. 

 Agent $i$ thus solves  the following strongly convex  problem:\vspace{-0.2cm}
 \begin{equation}
 \widetilde{\mathbf{x}}_{(i)}^n= \arg\min_{\mathbf{x}_{(i)}\in \mathcal{K}} \widetilde{F}_i\left(\mathbf{x}_{(i)};\mathbf{x}_{(i)}^n\right) + G^+\left(\mathbf{x}_{(i)}\right), \label{eq: x_tilde}\vspace{-0.2cm}
 \end{equation}
and updates its own local estimate $\mathbf{x}_{(i)}^n$  moving along the direction  $\widetilde{\mathbf{x}}_{(i)}^n- \mathbf{x}_{(i)}^n$ by a quantity (step-size) $\alpha^n>0$:\vspace{-0.2cm}
 \begin{equation}\label{wi}
 \mathbf{v}_{(i)}^n =\mathbf{x}_{(i)}^n + \alpha^n\left(\widetilde{\mathbf{x}}_{(i)}^n- \mathbf{x}_{(i)}^n\right). \vspace{-.3cm}
 \end{equation}
 
\vspace{.2cm}

\noindent\textbf{Step 2: Information mixing.}   We need   to introduce now
a mechanism to ensure that the local estimates  $\mathbf{x}_{(i)}^n$'s  eventually agree while $\widetilde{\boldsymbol{\pi}}_i^n$'s track  the gradients  $\sum_{j\neq i}\nabla f_j(\mathbf{x}_{(i)})$. 
Building on \cite{sun2016distributed}, consensus over   time-varying
digraphs without requiring the knowledge of the sequence of digraphs
and a double-stochastic weight matrix can be achieved employing the following broadcasting 
  protocol: given  $\mathbf{v}_{(i)}^n$, each agent $i$  updates its own
local estimate $\mathbf{x}_{(i)}^n$  together with one extra scalar variable  $\phi_i^n$ according to\vspace{-0.1cm} 
\begin{equation}
\phi_i^{n+1} =\!\sum_{j\in\mathcal{N}_i^{\textrm{in}}\left[n\right]}a_{ij}^n\phi_j^n,\quad\text{and}\quad
\mathbf{x}_{(i)}^{n+1} =\dfrac{1}{\phi_i^{n+1}}\sum_{j\in\mathcal{N}_i^{\textrm{in}}\left[n\right]}a_{ij}^n\phi_j^n\mathbf{v}_{(j)}^n,
\label{eq: update x} \vspace{-0.1cm}
\end{equation}
where $\phi_i^0=1$ for all $i$ and  $a_{ij}^n$'s are some weighting coefficients matching the graph $\mathcal{G}[n]$ in the following sense.\smallskip

\noindent\textbf{Assumption D (On the weighting matrix).} For all $n\geq 0$, the matrices $\mathbf{A}[n]\triangleq (a_{ij}^n)_{i,j}$ are chosen so that 
	\begin{enumerate*}
		\item[\textbf{(D1)}]  $a_{ii}^n\geq \kappa>0$ for all $i=1,\ldots,I$;
		\item[\textbf{(D2)}] $a_{ij}^n\geq \kappa>0$, if $\left(j,i\right)\in \mathcal{E}[n]$; and $a_{ij}^n=0$ otherwise;
		\item[\textbf{(D3)}] $\mathbf{A}[n]$ is column stochastic, i.e., $\mathbf{1} ^T\mathbf{A}[n]= \mathbf{1}^T$.
	\end{enumerate*} 
Some practical rules  satisfying the above assumption are given in Sec.\,\ref{subsec:discussion}. Here, we only remark  that, differently from most of the papers in the literature \cite{nedic2010constrained,nedic2009distributed,Lorenzo2016NEXT}, $\mathbf{A}[n]$ need not be double-stochastic but just column stochastic, which is a much weaker condition.\,This can be achieved thanks to the extra variables $\phi_i^n$ in \eqref{eq: update x}, whose goal roughly speaking is to dynamically build the missing    row-stochasticity, so that asymptotic consensus among  $\mathbf{x}_{(i)}^n$'s can be achieved. 

A similar scheme can be put forth  to update $\widetilde{\boldsymbol{\pi}}_i$'s building on the gradient tracking mechanism, first introduced in our work  \cite{Lorenzo2016NEXT}, and leveraging the information mixing protocol \eqref{eq: update x}, the desired update reads   (we omit further details because of the space limitation): each agent $i$ maintains  an extra (vector) variable  $\mathbf{y}_{(i)}^n$ [initialized as  $\mathbf{y}_{(i)}^n= \nabla f_i\left(\mathbf{x}_{(i)}^0\right)$];  and  consequently  $\widetilde{\boldsymbol{\pi}}_i$ are updated  according to \vspace{-0.2cm}
\begin{equation}
	\begin{aligned}
	&\mathbf{y}_{(i)}^{n+1}= \dfrac{1}{\phi_i^{n+1}} \sum_{j\in\mathcal{N}_i^{\textrm{in}}\left[n\right]}a_{ij}^n\left(\phi_j^n  \,\mathbf{y}_{(j)}^n + \nabla f_j(\mathbf{x}_{(j)}^{n+1}) - \nabla f_j(\mathbf{x}_{(j)}^{n})\right), \\
	& \widetilde{\boldsymbol{\pi}}_i^{n+1} = I \cdot \mathbf{y}_{(i)}^{n+1} - \nabla f_i\left(\mathbf{x}_{(i)}^{n+1}\right).\label{eq: update y}
	 \end{aligned}
	 \end{equation}
Note that the update of $\mathbf{y}_{(i)}$  and $\widetilde{\boldsymbol{\pi}}_i$  
 can be   performed locally by agent $i$, with the same
signaling as for \eqref{eq: update x}. 	
	One can show that if $\mathbf{x}_{(i)}^n$'s and $\mathbf{y}_{(i)}^n$'s are consensual (a fact that is proved in Th.\,1),  $\|\widetilde{\boldsymbol{\pi}}_i^n- \sum_{j\neq i}\nabla f_j\left(\mathbf{x}_{(i)}^n\right)\|\underset{n\rightarrow\infty}{\longrightarrow}0$. \\\indent We can now formally introduce the proposed algorithm, as given in Algorithm 1; its convergence properties are stated in Theorem 1 (the proof is omitted because of space limitation; see \cite{UnboundedGradient}). \vspace{-0.2cm}
	\begin{theorem}
	Let  $\left\{(\mathbf{x}_{(i)}^n)_{i=1}^{I}\right\}_n$  be the sequence generated by Algorithm 1, and let $\{\bar{\mathbf{z}}^n\triangleq \left(1/I\right)\sum_i \phi_{i}^n\mathbf{x}_{(i)}^n\}_n$.  Suppose that i) Assumptions A-D hold; ii) the step-size sequence $\{\alpha^n\}_n$ is chosen so that  $\alpha^n\in \left(0,1\right]$,  $\sum_{n=0}^{\infty} \alpha^n=+\infty$,	and  $\sum_{n=0}^{\infty} (\alpha^n)^2<+\infty$. Then, the following hold:	 (1)  $\bar{\mathbf{z}}^n$ is bounded for all $n$, and every limit point of $\bar{\mathbf{z}}^n$ is a d-stationary  solution of Problem \eqref{eq: P}; and  (2)  $\norm{\mathbf{x}_{(i)}^n-\bar{\mathbf{z}}^n}\to 0$ as $n\to +\infty$ for all $i$. \label{thm: convergence}\vspace{-0.4cm}
\end{theorem}
\begin{algorithm}
	\SetAlgoLined
	\KwData{For all agent $i$, $\mathbf{x}_{(i)}^0\in \mathcal{K}$, $\phi_i^0=1$, $\mathbf{y}_{(i)}^0=\nabla f_i\left(\mathbf{x}_{(i)}^0\right)$, $\widetilde{\boldsymbol{\pi}}_i^0= I\mathbf{y}_{(i)}^0-\nabla f_i\left(\mathbf{x}_{(i)}^0\right)$. Set $n=0$.}
	
	\CommentSty{[S.1]} If $\mathbf{x}_{(i)}^n$	satisfies termination criterion: STOP;\\
	\CommentSty{[S.2] Distributed Local SCA:}		
	Each agent $i$: \\
	\Indp (a) computes locally
	$\widetilde{\mathbf{x}}_{(i)}^n$ [cf.\,\eqref{eq: x_tilde}].
\\
	(b) updates its local variable $\mathbf{v}_{(i)}$ according to \eqref{wi}.
	
	\Indm \CommentSty{[S.3] Consensus: }
	Each agent $i$ broadcasts its local variables and sums up the received variables:\\
	\Indp 	 (a) Update $\phi_i^{n+1}$ and $\mathbf{x}_{(i)}^{n+1}$ using \eqref{eq: update x}.\\
	(b) Update   $\mathbf{y}_{(i)}^{n+1}$ and $\widetilde{\boldsymbol{\pi}}_i^{n+1}$ using  \eqref{eq: update y}.\\
	\Indm \CommentSty{[S.4] } $n\longleftarrow n+1$, go to \CommentSty{[S.1] }
	\caption{Distributed Sparse learning Algorithm (DSparsA)}\label{alg: DSL} 
\end{algorithm} 
 
Roughly speaking, Th.\,1 states two results:\,1)\,the weighted average $\bar{\mathbf{z}}^n$ of the $\mathbf{x}_i$'s converges to a d-stationary solution of  \eqref{eq: P}; 2) the $\mathbf{x}_i$'s  asymptotically agree on the common value $\bar{\mathbf{z}}^n$.  We remark that convergence is proved without requiring that  the (sub)gradients of $F$ or $G$ be bounded; this is a major achievement with respect to current distributed methods for nonconvex problems \cite{sun2016distributed,nedic2014stochastic,
nedic2015distributed,tatarenko2015non}.\vspace{-0.2cm}

\subsection{Discussion}\label{subsec:discussion}\vspace{-0.2cm}
\noindent\texttt{On the choice of $\tilde{f}_i$:}
 Assumption C is mild and offers a lot of flexibility in the choice of   $\tilde{f}_i$. Some examples are  the following:  

\noindent \textit{$-$Linearization:} One can always linearize $f_{i}$, which leads to  
  $\widetilde{f}_{i}(\mathbf{x}_{(i)};$ $\mathbf{x}_{(i)}^n) = f_{i}\left(\mathbf{x}_{(i)}^n\right)+  \nabla f_{i}\left(\mathbf{x}_{(i)}^n\right)^{\top}\left(\mathbf{x}_{(i)}-\mathbf{x}_{(i)}^n\right) +\frac{\tau_{i}}{2}\|\mathbf{x}_{(i)}-\mathbf{x}_{(i)}^n\|^{2}$.  

\noindent \textit{$-$Partial Linearization:} Consider the case that $f_{i}$ can be decomposed as $f_{i}(\mathbf{x}_{(i)})=f_{i}^{(1)}(\mathbf{x}_{(i)})+f_{i}^{(2)}(\mathbf{x}_{(i)})$, where $f_{i}^{(1)}$ is convex and $f_{i}^{(2)}$ is  nonconvex with Lipschitz continuous gradient.  Preserving the convex part of $f_i$ while linearizing $f_{i}^{(2)}$  leads to the following valid surrogate  
$
\widetilde{f}_{i}(\mathbf{x}_{(i)};\mathbf{x}_{(i)}^n) = f_{i}^{(1)}(\mathbf{x}_{(i)}) + f_{i}^{(2)}(\mathbf{x}_{(i)}^n)+ \frac{\tau_{i}}{2}\|\mathbf{x}_{i}-\mathbf{x}_{(i)}^n\|^{2} 
 + \nabla f_{i}^{(2)}(\mathbf{x}_{(i)}^n)^{\top}(\mathbf{x}_{i}-\mathbf{x}_{(i)}^n).$ We refer the readers to \cite{facchinei2015parallel,Lorenzo2016NEXT} for more choices of surrogates.
 \smallskip
%

\noindent\texttt{On the choice of the step-size.}  Several options are possible  for the  step-size sequence $\{\alpha^n\}_{n}$ satisfying the standard diminishing-rule in Th.\,1; see, e.g., \cite{bertsekas1999nonlinear}. Two instances  we found to be effective in our experiments are:  (1) $\alpha^n = \alpha_{0}/\left(n+1\right)^{\beta},$ with $\alpha_0>0$ and $0.5<\beta\leq 1$; and (2) $\alpha^n = \alpha^{n-1}\left(1-\mu\alpha^{n-1}\right),$ with $\alpha^0\in(0,1]$, and  $\mu\in\left(0,1\right)$.\smallskip

\noindent\texttt{On the choice of matrix $\mathbf{A}[n]$.}   In a digraph satisfying Assumption B, $\mathbf{A}[n]$ can be set to\vspace{-.2cm}
	\begin{equation}
	\begin{aligned}\label{eq: push sum weighting}
	a_{ij}^n = \begin{cases}
	1/d_j[n] & \left(j,i\right)\in \mathcal{E}[n],\\
	0 & \textrm{otherwise};
	\end{cases}
	\end{aligned}\vspace{-.1cm}
	\end{equation}
where $d_i[n] $ is the out-degree of agent $i$. Note that the message passing protocol in \eqref{eq: update x} and \eqref{eq: update y} based on (\ref{eq: push sum weighting})  can be easily implemented:  all agents only
need to i) broadcast their local variables 
normalized by their current out-degree; and ii) collect locally
the information coming from their neighbors.
Note that in the  special case that the \emph{undirected} graph,  $\mathbf{A}[n]$ becomes symmetric; consequently $\phi_i[n] = 1$ for all $i=1,\ldots,I$ and $n\in\mathbb{N}_+$ [i.e., the update of $\phi$ in step \eqref{eq: update x}  can be eliminated], and Assumption D3 is readily satisfied choosing $\mathbf{A}[n]$ according to rule proposed in the literature for    double-stochastic matrices; some widely used rules are: the uniform weights \cite{blondel2005convergence}, Laplacian weights \cite{scherber2004locally}, and the Metropolis-Hastings weights \cite{xiao2005scheme}.\vspace{-0.2cm}

\section{Numerical Results}\vspace{-0.2cm}
In this section we test DSparsA on two instances of Problem \eqref{eq: P}, namely: i) the sparse linear regression problem \eqref{p: sparse regression} with the ``Log'' penalty function; and ii)  the sparse PCA problem \eqref{p: sPCA} with SCAD penalty function given in Table\,\ref{tab: ncvx regularizer}. For both problems, we simulated a network composed of $I=30$ users; the sequence of time-varying digraphs is such that, at   each time slot the graph is strongly connected and  every agent has two out-neighbors.\smallskip

\noindent \textbf{Example\,\#1:\,Sparse regression.} Consider Problem\,\eqref{p: sparse regression}, where $G$ is the ``Log'' penalty function given in Table\,\ref{tab: ncvx regularizer}.  The underlying sparse linear  model is  
 $\mathbf{b}_i = \mathbf{A}_i\mathbf{x}_0 +\mathbf{n}_i$, where $\mathbf{A}_i\in \mathbb{R}^{20\times 200}$  is the sensing matrix (with rows normalized to one),  $\mathbf{x}_0\in \mathbb{R}^{200}$ is  the unknown signal, and $\mathbf{n}_i$ is the observation noise, all randomly generated, with  i.i.d Gaussian entries. Each component of the noise vector has  standard deviation $\sigma_i=0.1$. To impose sparsity on   $\mathbf{x}_0$, we set, uniformly at random, to zero  $80\%$ of  its component.   Finally, we set  $\theta = 20$ [cf.\,Table\,1] and  $\lambda = 0.5$. We tested the following two instances  of DSparsA: i)  DSparsA-SCA,  wherein the surrogate function $\widetilde{f}_i$ coincides with  $f_i$, since $f_i$ is already convex. In this case, to compute  $\widetilde{\mathbf{x}}_{(i)}^n$, each agent needs to solve a LASSO problem, which can be efficiently done using the FLEXA  algorithm \cite{facchinei2015parallel}; and ii)   DSparsA-L,  where $\widetilde{f}_i$ is constructed  by linearizing   $\widetilde{f}_i$ at  $\mathbf{x}_{(i)}^n$, as shown in Sec.\,\ref{subsec:discussion}, where   %
 $\nabla f_i\left(\mathbf{x}_{(i)}^n\right) = 2\mathbf{A}_i^\top\left(\mathbf{A}_i\mathbf{x}_{(i)}^n - \mathbf{b}_i\right)$. Consequently, the update $\widetilde{\mathbf{x}}_{(i)}^n$ has the following   closed form \vspace{-0.2cm}
\begin{equation*}
	\widetilde{\mathbf{x}}_{(i)}^n = \mathcal{S}_{\frac{\eta\lambda}{\tau_i}}\left\{\mathbf{x}_{(i)}^n-\frac{1}{\tau_i}\left(\nabla f_i\left(\mathbf{x}_{(i)}^n\right) + \widetilde{\boldsymbol{\pi}}_i^n - \lambda\cdot\nabla G^-\left(\mathbf{x}_{(i)}^n\right)\right)\right\},\vspace{-0.1cm}
\end{equation*}
where $\mathcal{S}_{\eta\lambda/\tau_i}(\bullet)$ is the soft-thresholding operator, {defined as $\mathcal{S}_\lambda\left(x\right) \triangleq \textrm{sign}\left(x\right)\cdot\max\left(\lvert x\rvert-\lambda,0\right)$, and the explicit expression of  $\eta(\theta)$ and $\nabla G^-$ is given in row 5 of Table\,\ref{tab: gradient}}.

Since there are no convergent distributed schemes in the literature for the problem under consideration,  we compare our algorithms  with the subgradient-push algorithm \cite{nedic2015distributed}, developed for convex functions with bounded subgradients. We report results achieved with the following  tuning of the algorithms (which provided the best performance, among all the choices we tested).  For all algorithms, we used  the step-size rule  (2), as given in Sec.\,\ref{subsec:discussion}: in DSparsA, we set $\alpha^0 = 0.1$ and $\mu = 10^{-3}$ whereas for the subgradient-push we used $\alpha^0=1$ and $\mu = 10^{-2}$. Finally, in DSparsA, we set $\tau_i=2$ for all $i$. We use two merit functions to measure the progresses of the algorithms, namely: i)   $J^n \triangleq ||\bar{\mathbf{z}}^n-\mathcal{S}_{\eta\lambda}\{\bar{\mathbf{z}}^n-(\nabla F(\bar{\mathbf{z}}^n) - \lambda\cdot\nabla G^-(\bar{\mathbf{z}}^n))\}||_\infty$, which measures the distance of the weighted average $\overline{\mathbf{z}}^n$ from d-stationarity (note that such a function is zero if and only if the argument is a d-stationary solution of \eqref{p: sparse regression}); and ii)$D^n\triangleq \max_{i}\{\|\mathbf{x}_{(i)}^n - \bar{\mathbf{z}}^n\|_\infty\}$ , which measures how far are the agents to reach consensus. 
In Fig. 1(a) [resp.\,Fig.\,1(b)] we plot $D^n$ and $J^n$ [resp. the normalized MSE, defined as   $\textrm{NMSE}^n = \frac{1}{I}\sum_{i=1}^{I}||\mathbf{x}_{(i)}^n - \mathbf{x}_0||_2^2/||\mathbf{x}_0||_2^2$] achieved by all the algorithms vs. the  total number of
communication exchanges per node. For the subgradient-push algorithm, this
number coincides with the iteration index $n$ whereas  for DSparsAs it  is $2n$.   All the curves are averaged over $100$ independent noise
realizations. The
figures clearly show that both versions of DSparsA are much
faster than the subgradient-push algorithm (or, equivalently, they require
less information exchanges), which is not even guaranteed to converge.      Moreover, as expected,  DSparsA-SCA reaches
high precision faster than DSparsA-L; this is mainly due to
the fact that the surrogate function in the former retains the partial
convexity of  $f_i$ rather than just linearizing it.

\begin{figure}[t]
\vspace{-0.7cm}	\centering
	\begin{subfigure}[t]{0.21\textwidth}
		\includegraphics[scale = .14]{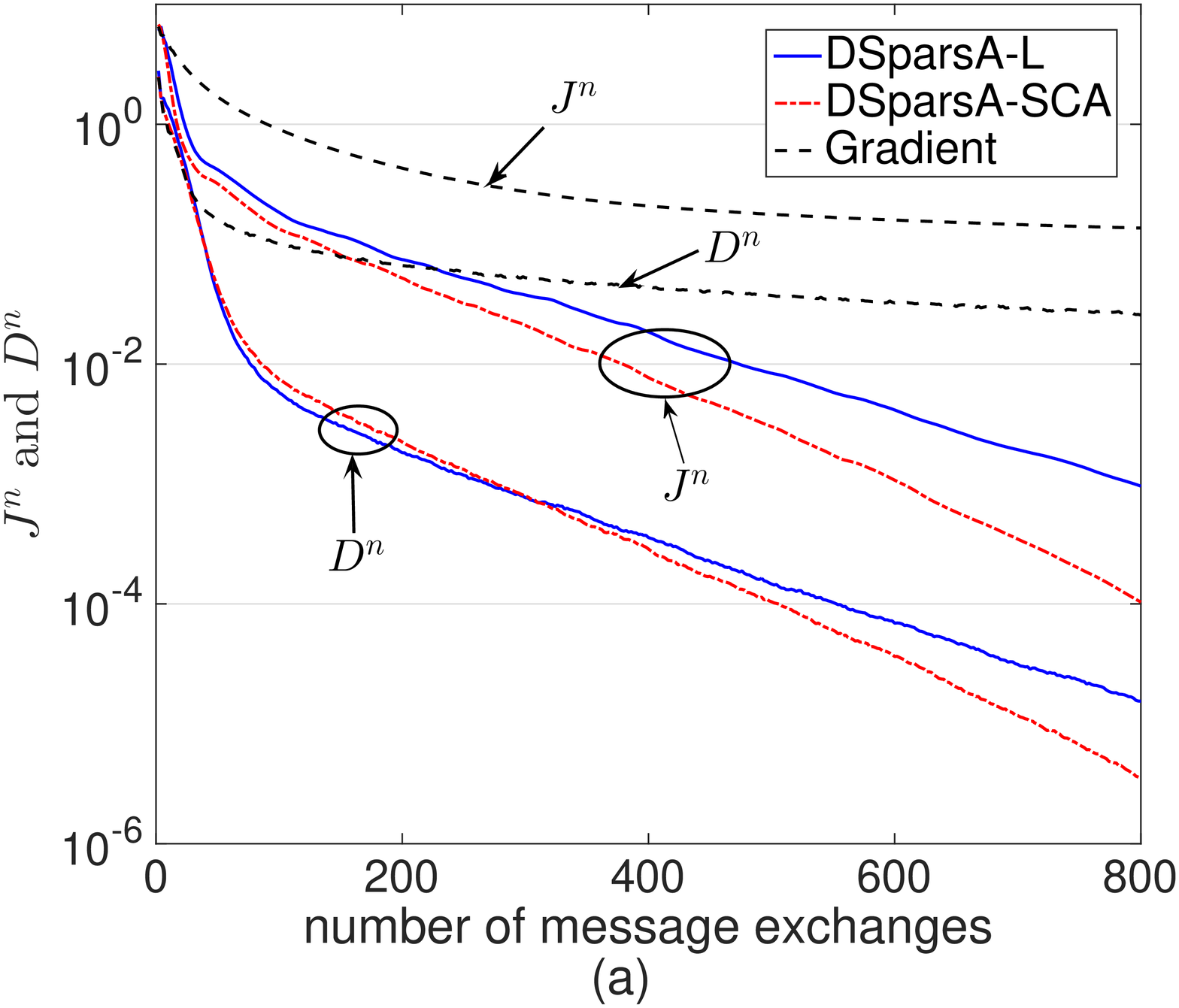}

		\label{fig: opt criteria-log}	
	\end{subfigure}
	~
	\begin{subfigure}[t]{0.21\textwidth}
		\includegraphics[scale=.15]{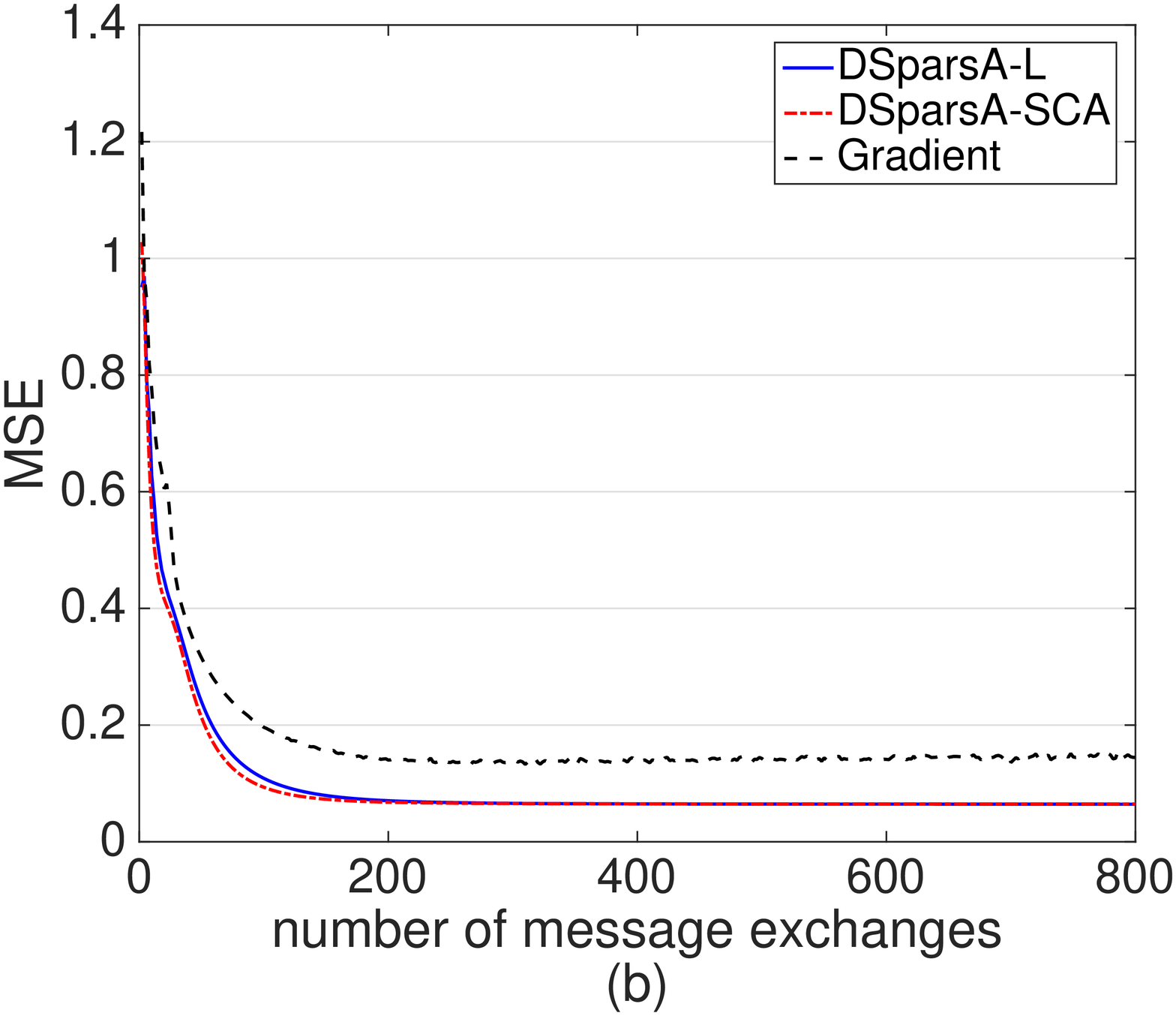}
\label{fig: MSE-log}
	\end{subfigure}\vspace{-0.7cm}
	\caption{\small Sparse linear regression:   Optimality  measurements  $J^n$ and  consensus disagreement $D^n$ [subplot (a)] and   normalized MSE $\textrm{NMSE}^n$ [subplot (b)] versus  per-node communication exchanges.}\vspace{-0.3cm}
	\label{fig: log-penalty}
\end{figure}

\noindent \textbf{Example\#2: Distributed Sparse PCA.} Consider the sparse PCA problem \eqref{p: sPCA}, where $G$ is the SCAD penalty function [cf.\,Table\,1];    $\theta = 20$, $a= 2$, and  $\lambda = 5$.  The rows of  data matrix $\mathbf{D}_i\in \mathbb{R}^{500\times30}$ are generated as  i.i.d  Gaussian random vectors of mean zero and covariance $\boldsymbol{\Sigma}$. 
 The leading eigenvector $\mathbf{u}_1$of $\boldsymbol{\Sigma}$ with eigenvalue $12$ is dense, while the next two eigenvectors $\mathbf{u}_2$ and $\mathbf{u}_3$ are of cardinality $5$ with eigenvalues being $10$ and $8$, respectively. The rest of the $\mathbf{u}_i$'s are randomly generated with eigenvalue less than $5$. The task is to estimate $\mathbf{u}_2$ from the $\mathbf{D}_i$'s.
Since $f_i $ is concave, the surrogate function $\tilde{f}_i$ is obtained by linearizing  $f_i $ [cf.\,Sec.\,\ref{subsec:discussion}]. The convexified optimization problem of agent $i$'s   reads\vspace{-0.2cm}
 \begin{equation}
 \begin{aligned}
 &\underset{\norm{\mathbf{x}}_2 \leq 1}{\textrm{min}}& & \mathbf{g}_i^{n\top}\mathbf{x} + \frac{\tau_i}{2}\norm{\mathbf{x} -\mathbf{x}_{(i)}^n}^2 + \lambda G^+\left(\mathbf{x}\right),\\
 \end{aligned} \label{p: bisection}\vspace{-0.2cm}
 \end{equation}
 where $\mathbf{g}_i^n= \nabla f_i\left(\mathbf{x}_{(i)}^n\right) + \widetilde{\boldsymbol{\pi}}_i^n - \nabla G^-\left(\mathbf{x}_{(i)}^n\right)$, and we set $\tau_i = 10^{-3}$. The unique solution $\widetilde{\mathbf{x}}^n_i$of (\ref{p: bisection}) can be efficiently obtained using the soft-thresholding operator, followed by a scalar bi-section; we omit the details because of space limitation. 
  
We compare DSparsA with subgradient-push, where we added a projection step after the gradient descent to maintain feasibility. Note that there is no formal proof of convergence for such an algorithm.   For both algorithms, we used  the step-size rule  (2), as given in Sec.\,\ref{subsec:discussion}: in DSparsA, we set   $\alpha^0 = 1$, and $\mu = 10^{-3}$ whereas for the subgradient-push we used $\alpha^0=0.1$ and $\mu = 10^{-2}$. Since Problem \eqref{p: sPCA} is constrained, we modify the stationarity measure as $J^n \triangleq \|\widehat{\mathbf{x}}(\bar{\mathbf{z}}^n)-\bar{\mathbf{z}}^n\|_{\infty}$, with
$\widehat{\mathbf{x}}(\bar{\mathbf{z}}^n)\triangleq 
  {\text{argmin}}_{\|\mathbf{x}\|\leq1} \{\lambda G^+\left(\mathbf{x}\right)+ \left(\nabla F\left(\bar{\mathbf{z}}^n\right)-\lambda\nabla G^-\left(\bar{\mathbf{z}}^n\right)\right)^\top \mathbf{x} + \norm{\mathbf{x}-\bar{\mathbf{z}}^n}^2/2\}.$ In Fig.\,2(a) [resp. Fig.\,2(b)] we plot $D^n$ and $J^n$ [resp.\,$\textrm{NMSE}^n$, which is defined as in Example \#1, with $\mathbf{x}_0$ replaced by $\mathbf{u}_2$] achieved by all the algorithms versus the  total number of
communication exchanges per node.\,All the curves are averaged over $100$ independent data generations. The figures show that DSparsA significantly outperforms the subgradient   method both in terms of convergence speed  and MSE. 
\vspace{-0.3cm}
 \begin{figure}
 
 	\begin{subfigure}[t]{0.22\textwidth}
 		\includegraphics[scale = .20]{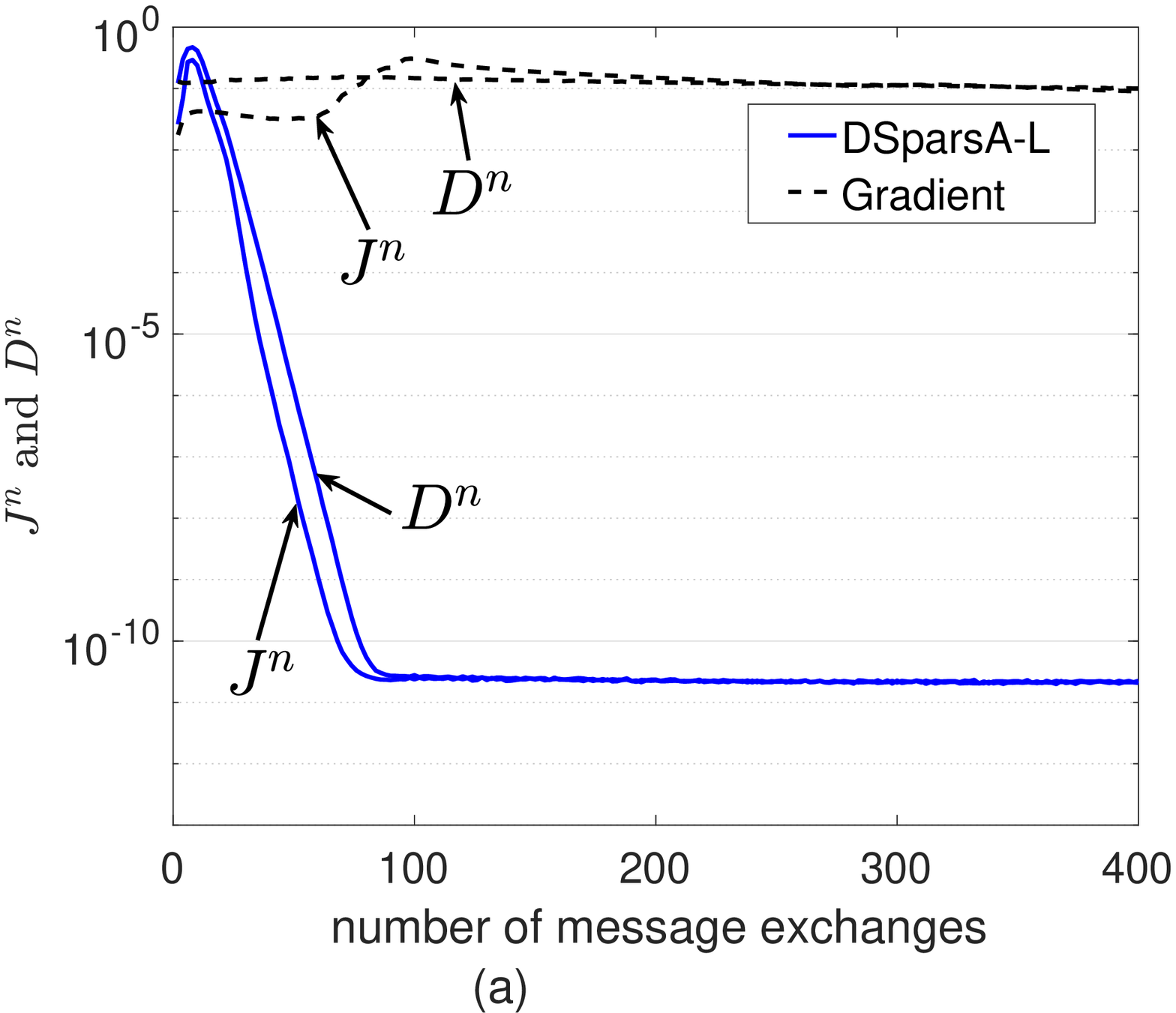}
 		
 		\label{fig: opt criteria-log}	
 	\end{subfigure}
 	~
 	\begin{subfigure}[t]{0.21\textwidth}
 		\includegraphics[scale=.2]{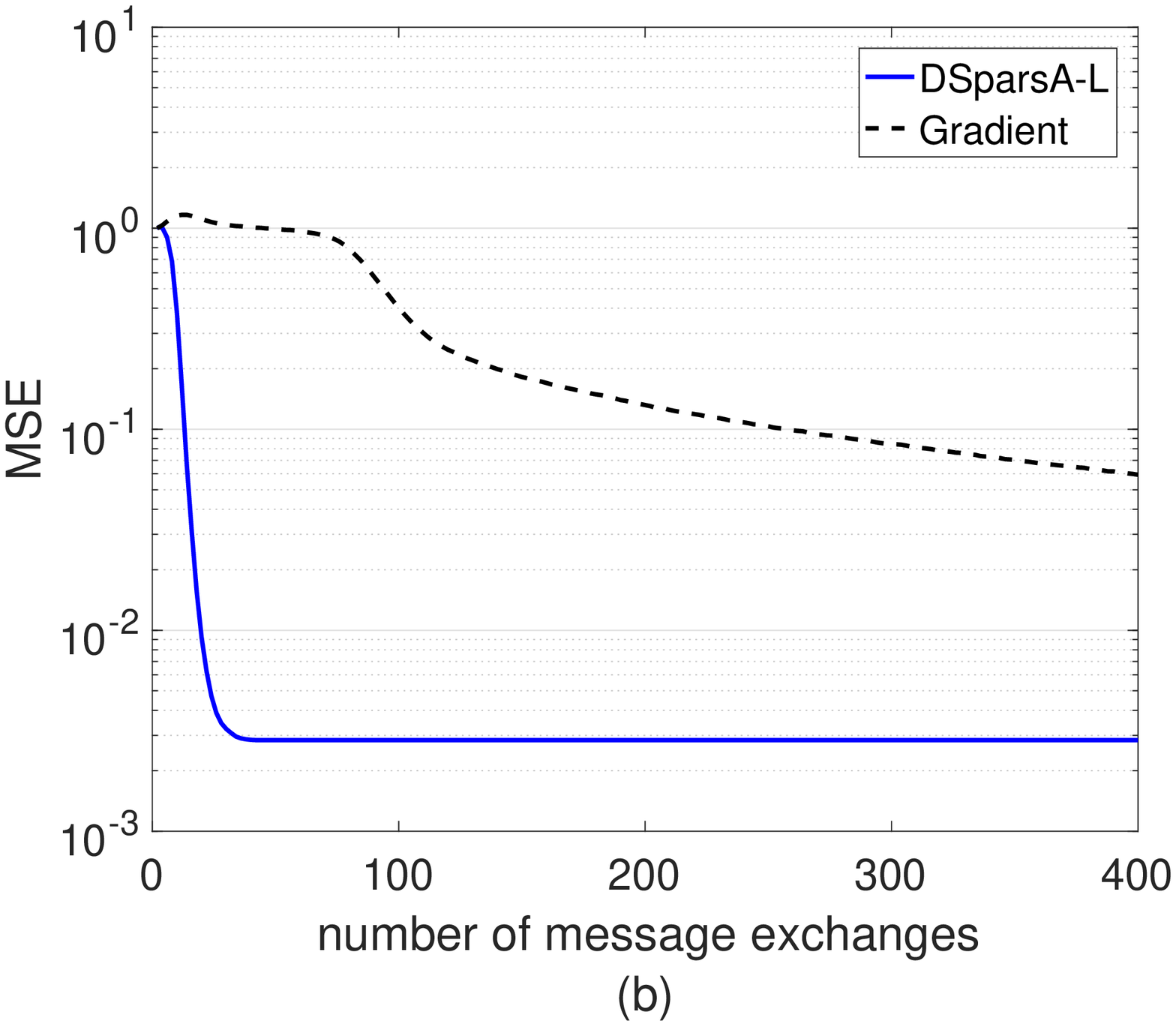}
 		\label{fig: MSE-log}
 	\end{subfigure}\vspace{-0.6cm}
 	\caption{\small Sparse PCA:   Optimality  measurements  $J^n$ and  consensus disagreement $D^n$ [subplot (a)] and   normalized MSE $\textrm{NMSE}^n$ [subplot (b)] versus  per-node communication exchanges.}\vspace{-0.5cm}
 	\label{fig: log-penalty}
 \end{figure}
 \section{conclusions}\vspace{-0.2cm}
 We have proposed the first unified  distributed  algorithmic framework for the computation     of d-stationary solutions of a fairly general class of non convex statistical learning problems. Our scheme is implementable over time-varying network with arbitrary topology and does not require that the (sub)gradient of the objective function is bounded.

\bibliographystyle{IEEEbib}
\footnotesize
\bibliography{refs.bib}

\end{document}